\documentclass[a4paper,10pt]{article}
\usepackage[a4paper, margin=1.2in]{geometry}
\usepackage{mathrsfs,amsfonts,latexsym,theorem,amssymb,amsmath,newlfont}
\usepackage{graphicx}
\usepackage{dsfont}
\usepackage{pgf}
\usepackage{srcltx}
\usepackage{fancyhdr}
\usepackage{times}
\usepackage{pstricks,pstricks-add,pst-math,pst-xkey} 
\usepackage{nicefrac}
\usepackage{mathtools}
\usepackage{hyperref}

\DeclarePairedDelimiter{\ceil}{\lceil}{\rceil}

\theoremstyle{plain}

\newtheorem{thm}{Theorem}[section]
\newtheorem{defn}[thm]{Definition}
\newtheorem{lem}[thm]{Lemma}
\newtheorem{prop}[thm]{Proposition}
\newtheorem{cor}[thm]{Corollary}
\newtheorem{remark}[thm]{Remark}
\newtheorem{open}[thm]{Open Problem}

\def\beginpf{\noindent {\bf Proof:} \quad}
\def\endpf{\rightline{$\square$}}

\def\NN{{\mathbb N}}

\def\RR{{\mathbb R}}

\def\N{{\cal N}}

\def\O{{\cal O}}

\def\LL{{\cal L}}

\def\<{\langle}
\def\>{\rangle}

\DeclarePairedDelimiter{\floor}{\lfloor}{\rfloor}
\DeclareMathOperator{\sech}{sech}

\title{Lagrangian Controllability of the 1-D Korteweg-de Vries Equation}
\author{Ludovick \textsc{Gagnon}\footnote{Universit\'e Pierre et Marie Curie, LJLL UMR 7598, 4 place Jussieu, 75005, Paris, France.}}

\begin{document}

\maketitle

\begin{abstract}

We consider in this paper the problem of the Lagrangian controllability for the Korteweg-de Vries equation. Using the $N$-solitons solution, we prove that, for any length of the spatial domain $L>0$ and any time $T>0$, it is possible to choose appropriate boundary controls of KdV equation such that the flow associated to this solution exit the domain in time $T$. 

\end{abstract}

\section{Introduction}

The Korteweg-de Vries (KdV) equation defined on the real line 
\begin{equation}\label{KdVReal}
y_t+y_x+y_{xxx}+yy_x=0, \quad x\in \RR, t\in \RR,
\end{equation}
first derived, independently, by Boussinesq in 1877 (\cite{Boussinesq}) and by Korteweg and de Vries in 1895 (\cite{KortOriginal}), is obtained as a first order approximation of the free-surface solution of the full governing equations for a homogeneous, non-viscous and irrotational shallow fluid (see, for example, \cite[p. 460]{Whit} for a complete derivation of the KdV from the governing equations). In the context of water waves, solutions of (\ref{KdVReal}) correspond to the free-surface of approximately two-dimensional waves (the motion of the waves are assumed to be parallel to the crest) travelling from left to right. More recently, the KdV equation has found applications in the context of collisionless plasma hydromagnetic waves \cite{Gardner}, long waves in anharmonic crystals \cite{Zabusky}, ion-acoustic plasma \cite{Washimi} and cosmology \cite{Lidsey}. 

In this paper, we are interested by the small-time Lagrangian controllability of the Korteweg-de Vries equation starting from rest
\begin{align}
y_t+y_x+y_{xxx}+yy_x&=0 ,  &x&\in [0,L], \quad t\in [0,T], \label{KdV} \\ 
y(0,t)&=u(t) ,  &t&\in [0,T],  \label{KdVu} \\
y(L,t)&=v(t) ,  &t&\in [0,T],  \label{KdVv} \\
y_x(L,t)&=w(t) ,  &t&\in [0,T], \label{KdVw} \\
y(x,0)&=0, &x&\in [0,L], \label{KdV0} 
\end{align}
with boundary controls $u(t), v(t)$ and $w(t) \in \RR$ and $T,L>0$. The Lagrangian controllability of (\ref{KdV})-(\ref{KdV0}) is defined as follow,

\begin{defn}[Small-time Lagrangian Controllability]\label{LagDef}
Equations (\ref{KdV})-(\ref{KdV0}) are small-time Lagrangian controllable if and only if, for all $T,L>0$, there exists $u(t),v(t)$ and $w(t)\in \RR$ such that, if we consider $\hat{y}$ the extension of the solution $y$ of (\ref{KdV})-(\ref{KdV0}) by 
\begin{equation}\label{extensionlag}
 \hat{y}(x,t)=\left\{
 \begin{array}{l l}
y(0,t), & \textrm{ if } x\leq 0,\\ 
y(x,t), & \textrm{ if } x\in [0,L], \\
y(L,t), & \textrm{ if } x\geq L,
\end{array}
\right.
\end{equation}
the flow $\Phi$ defined by
\begin{equation}\label{flux}
\left\{
\begin{array}{rcll}
\dfrac{\partial \Phi}{\partial t}(x,t) &=&\hat{y}(\Phi(x,t),t),& x\in \RR, t\in \RR^+, \\
\Phi(x,0)&=&x, &x\in \RR,
\end{array}
\right.
\end{equation}
satisfies $\Phi(x,T)\geq L, \, \forall x\in [0,L]$.
\end{defn}

For sake of simplicity, we will always refer, in the following, to this construction when speaking of the flow of a solution of (\ref{KdV})-(\ref{KdV0}).

Let us describe a physical interpretation of the Lagrangian controllability in the water waves context. Consider the Cartesian coordinates $(x,z)$ such that $x$ is the horizontal direction in which the waves travel, $z=0$ denotes the flat bottom of the fluid, $z=h_0$ is the height of the fluid at rest and $z=h_0+y(x,t)$ is the free-surface where $y$ is solution of (\ref{KdV})-(\ref{KdV0}). The flow defined by (\ref{flux}) is, up to a physical constant, the horizontal component (independent of $z$) of first order approximation of the velocity field of (\ref{KdVReal}) (\cite[p. 460]{Whit}). Therefore, one may interpret the problem of the Lagrangian controllability of (\ref{KdV})-(\ref{KdV0}) as the problem of moving the particles initially located in the region $[0,L]\times [0,h_0]$ at time $t=0$ to the right of $L$ at time $t=T$ by means of waves created by the boundary controls. A possible application of the Lagrangian controllability is the displacement of polluted water in channels to a waste water treatment plant. 

\begin{remark}\label{remark}
It is a natural condition to impose that the flow exit to the right and not to the left since solutions of the KdV equation correspond to waves travelling from left to right. This assumption in the derivation of the KdV equation (\ref{KdVReal}) is transposed in the asymptotic behaviour, given by the Inverse Scattering Method (\cite{Segur}), of solutions of (\ref{KdVReal}) for smooth initial data : a finite number of solitons travelling to the right and a decaying wave train to the left. This asymmetric behaviour is different, for example, from the Euler equation for which there would be no geometrical restriction on where the flow should leave for a similar problem. In fact, the solution constructed for the return method by Coron and Glass to show the global Eulerian controllability of the 2-D and 3-D Euler equation respectively (\cite{CoronEuler},\cite{GlassEuler}) also provides a proof that the flow exit the domain for the 2-D and 3-D Euler equations.
\end{remark}

The main result of this paper is the small-time Lagrangian controllability of (\ref{KdV})-(\ref{KdV0}), with the additional property that the solution is at rest at time $t=T$. 
\begin{thm}\label{Lag}
Let $T,L>0$. Then, there exists $y\in C([0,T];H^2(0,L))$ solution of (\ref{KdV})-(\ref{KdV0}) such that the associated flow $\Phi$ satisfies 
\begin{align} 
\dfrac{\partial \Phi}{\partial t}(x,0)=& \, 0, \qquad x\in [0,L], \label{p1} \\
\dfrac{\partial \Phi}{\partial t}(x,T)=& \, 0, \qquad x\in [0,L], \label{p2} \\
\Phi(x,t)  \geq & \, L,  \qquad x\in [0,L], \, t\in [T,\infty). \label{p4} 
\end{align}
\end{thm}

Theorem \ref{Lag} follows from an explicit solution of (\ref{KdVReal}) satisfying (\ref{p4}) and a smallness condition on the state in the neighborhood of $t=0$ and $t=T$. 

\begin{thm}\label{intermediaire}
Let $L,T,\delta>0$ and $(\epsilon_1,\epsilon_2)\in (0,T/2)^2$. Then there exists a positive solution $y\in C^{\infty}([0,T]\times \RR)$ of (\ref{KdVReal}) such that
\[
\|y(.,t)\|_{H^2(0,L)}<\delta, \quad \forall t\in (0,\epsilon_1)\cup(T-\epsilon_2,T),
\] and such that the flow $\Phi(x,t)$ associated to $y$ 
satisfies 
\[\Phi(x,t)\geq L, \quad \forall (x,t)\in [0,L]\times (T-\epsilon_2,T).\]
\end{thm}

This solution is constructed by means of the $N$-solitons solution. For later conveniences, let us express the soliton solution of KdV for the change of variables $x \mapsto x-t$, $y \mapsto 6\eta$. The Korteweg-de Vries equation becomes 
\begin{equation}\label{KdVHirota}
\eta_t+6\eta \eta_x+\eta_{xxx}=0, \quad x\in \RR, t\in \RR,
\end{equation}
and solitons of (\ref{KdVHirota}) are given, for $\alpha>0$ and $s\in \RR$, by
\begin{equation}\label{soliton}
\eta(x,t)=\dfrac{\alpha^2}{2}\sech^2\left(\dfrac{-\alpha(x-s)+\alpha^3t}{2}\right). 
\end{equation}
We note the following for solitons given by (\ref{soliton})
\begin{enumerate}
\item The amplitude is given by $\alpha^2/2$ and is reached at $x=s+\alpha^2t$;
\item The travelling speed is $\alpha^2$;
\item The distance between the $x$ coordinates where the height of the soliton is $\alpha/4$,  defined as the width, is 
\begin{equation}\label{width}
w(\alpha):=\frac{4}{\alpha}\ln(\sqrt{2\alpha}(1+\sqrt{1-\frac{1}{2\alpha}})).
\end{equation}
\end{enumerate}
From these properties, one remarks that taller solitons travel faster and are narrower. Moreover, the width of a soliton tend to infinity as its amplitude tends to zero. 

We point out that Theorem \ref{Lag} is not the consequence of the passage of a single soliton inside the domain $[0,L]$. Indeed, consider the flow defined on the whole real line
\[
\dfrac{\partial \Phi}{\partial t}(x,t)=\eta(\Phi(x,t),t), \quad x\in \RR, t\in \RR,
\]
where $\eta$ is a soliton, with $\alpha>0$ and $s\in \RR$, solution of (\ref{KdVHirota}). Then, from (\ref{soliton}), one obtains that the total displacement $\left| \Phi(x,\infty)-\Phi(x,-\infty)\right|$ is of order $2/\alpha$ and, since the speed of propagation of a soliton is $\alpha^2$, one cannot obtain, at the same, time both a large displacement and an arbitrarily small time. 

We rather use the $N$-solitons solution to prove Theorem \ref{Lag}. Expressed in the closed form by Hirota in 1971 \cite{Hirota}, the $N$-solitons solution of (\ref{KdVHirota}) writes

\begin{align}\label{Nsolitons}
\eta&= -2\left( \ln F\right)_{xx}, \\ 
F&= 1+\sum_{n=1}^N\sum_{C_{n}^{N}}a(i_1,...,i_n) \prod_{j=1}^{n} f_{i_j}, \nonumber
\end{align}
where, for $1\leq i \leq N$, $N\in \NN$
\[
 f_i(x,t)=\exp\left(-\alpha_i(x-s_i)+\alpha_i^3t\right), \, s_i\in \RR, \alpha_i>0, 
\]
where
\begin{align*}
a(i_1,...,i_n)&= \prod^{n}_{k<l} a(i_k,i_l), \qquad   \textrm{ if } n\geq 2, \\
a(i_k,i_l)&=\left(\dfrac{\alpha_{i_k}-\alpha_{i_l}}{\alpha_{i_k}+\alpha_{i_l}}\right)^2, \\
a(i_k)&=1,  				
\end{align*}
where $\sum_{C_{n}^{N}}$ is the sum over all the $n$ indexes $i_1,...,i_n$, taken, without permutations, from $\{1,\ldots,N\}$. 

To explain why this solution is called the $N$-solitons solution, let us consider the case where $N=2$. The solution is expressed as
\[
\dfrac{\eta}{2}=\dfrac{\alpha_1^2 f_1 + \alpha_2^2 f_2 + 2(\alpha_2-\alpha_1)^2f_1f_2+((\alpha_2-\alpha_1)/(\alpha_2+\alpha_1))^2(\alpha_2^2 f_1^2f_2+ \alpha_1^2f_2^2 f_1)}{(1+f_1+f_2+((\alpha_2-\alpha_1)/(\alpha_2+\alpha_1))^2f_1f_2)^2}.
\]
If $f_1\simeq 1$, $f_2 \ll 1$, the behaviour of $\eta$ is then given by
\begin{align*}
\eta &\simeq  2\dfrac{\alpha_1^2 f_1}{(1+f_1)^2} \\
&= \dfrac{\alpha_1^2}{2}\sech^2\left(\dfrac{-\alpha_1(x-s_1)+\alpha_1^3t}{2}\right),
\end{align*}
while, in the case where $f_2\simeq 1$, $f_1 \gg 1$, we have, 
\begin{align*}
\eta &\simeq 2\dfrac{\alpha_2^2((\alpha_2-\alpha_1)/(\alpha_2+\alpha_1))^2 f_1^2f_2}{(f_1+((\alpha_2-\alpha_1)/(\alpha_2+\alpha_1))^2f_1f_2)^2} \\
&=  2\dfrac{\alpha_2^2\exp(\ln(((\alpha_2-\alpha_1)/(\alpha_2+\alpha_1))^2))f_2}{(1+\exp(\ln(((\alpha_2-\alpha_1)/(\alpha_2+\alpha_1))^2))f_2)^2} \\
&= \dfrac{\alpha_2^2}{2}\sech^2\left(\dfrac{-\alpha_2(x-s_2-\frac{1}{\alpha_2}\ln(((\alpha_2-\alpha_1)/(\alpha_2+\alpha_1))^2))+\alpha_2^3t}{2}\right),
\end{align*}
that is, a soliton with a phase shift of $\frac{1}{\alpha_2}\ln(((\alpha_2-\alpha_1)/(\alpha_2+\alpha_1))^2)$. 

Let us now describe the behaviour of the $2$-solitons solution in the case where $0<\alpha_1<\alpha_2$ and let us denote them soliton 1 and 2 respectively. When $t\rightarrow -\infty$, the solitons behaves like 2 distinct solitons, soliton 1 being ahead of soliton 2 and the latter having a phase shift of $\frac{1}{\alpha_2}\ln(((\alpha_2-\alpha_1)/(\alpha_2+\alpha_1))^2)$. When $f_1\simeq f_2 \simeq 1$, interactions occur between soliton 1 and 2 and, during this period, their combined amplitude decreases. Moreover, if $\alpha_1$ and $\alpha_2$ are of same magnitude, they exchange their amplitudes and velocities (\cite{Zabusky2}). After the interaction, they behave as distinct solitons and are left unchanged in shape, the only notable effect of the interaction is the phase shift of soliton 1, of $\frac{1}{\alpha_1}\ln(((\alpha_2-\alpha_1)/(\alpha_2+\alpha_1))^2)$, while the soliton 2 no longer has one. Those effects are the result of the nonlinearity of the equation. Figure 1 illustrates the phase shift produced when two solitons interact for the 2-solitons solution. The frame is fixed at the speed of the slower soliton. 

\begin{figure}[!th]
\begin{center}
	\includegraphics[height=5cm]{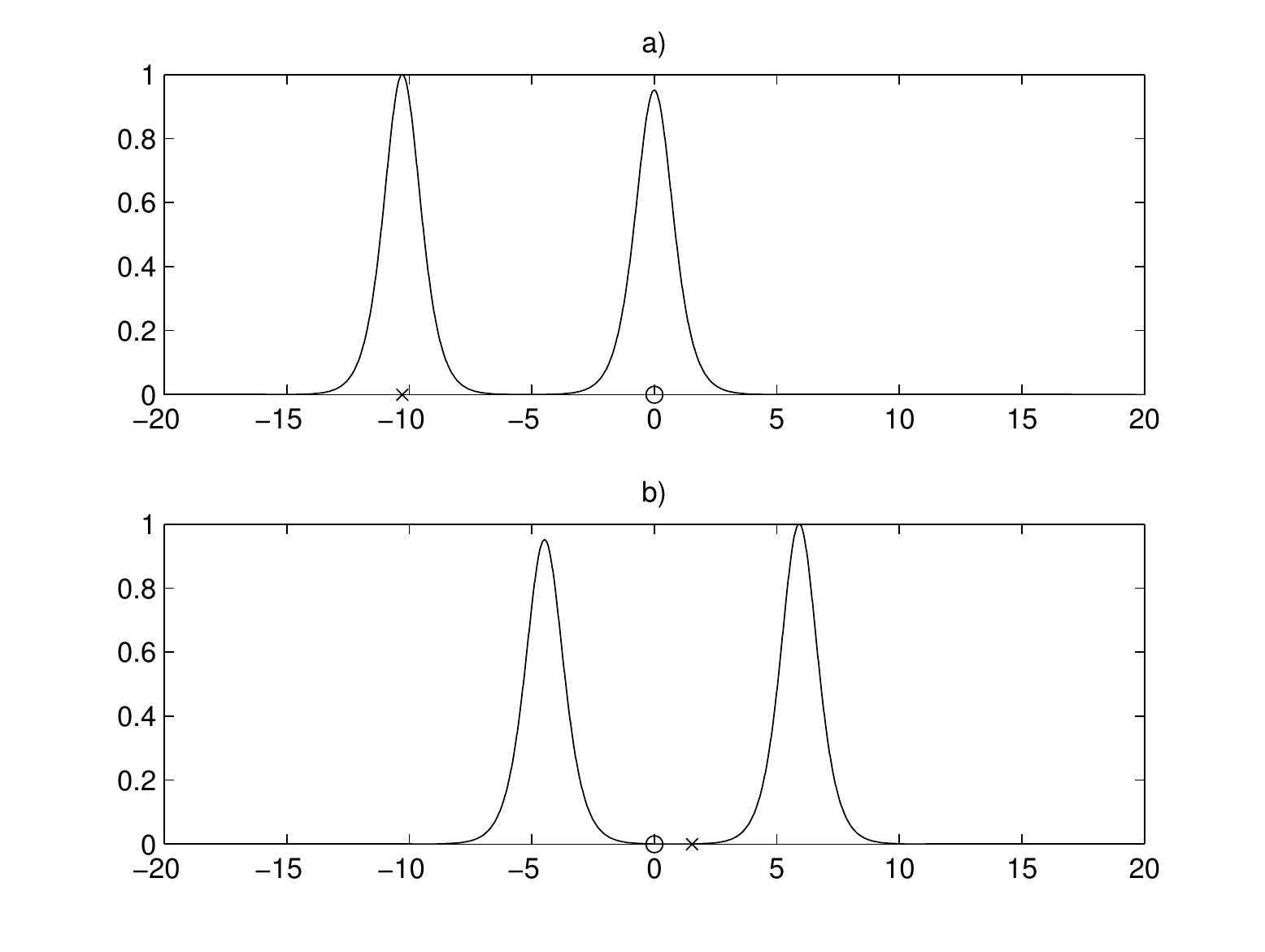}
\end{center}
\caption{\footnotesize{An interaction between two solitons. The cross (circle) represents the position of the maximum of faster (slower) soliton if no interaction would have occured. Figure a) is the state of the solution before the collison and b) is after the collision}}\label{interaction}
\end{figure}

Let us sketch the three steps of the proof of Theorem \ref{intermediaire}. Consider the $N$-solitons solution.

\begin{enumerate}
\item[Step 1.] At time $t=0$, the solitons are located at the left of $x=0$. They are ordered in increasing order of height to prevent further interactions. During this step, only the tail of the $N$-solitons solution is located inside the domain $[0,L]$. An estimation of the norm of the tail with respect of $\alpha_i$ and $N$ shows that the larger the $\alpha_i$ are, the smaller the norm of the tail is. 
\item[Step 2.] During the time interval $(0,T)$, the $N$ solitons travel inside the domain. A lower bound of the flow is provided with respect of $\alpha_i$ and $N$. This lower bound is estimated by the displacement induced by the $N$ solitons, the displacement due to the interactions being negligible. 
\item[Step 3.] At time $t=T$, the solitons are located at the right of $x=L$. Only the tail of the $N$-solitons solution is located in $[0,L]$. A similar estimate than in step 1 proves that the norm of the tail is small for $\alpha_i$ large. 
\end{enumerate}

With the required estimates at hand, one chooses $\alpha_i$ and $N$ large enough to obtain Theorem \ref{intermediaire}. 

To prove Theorem \ref{Lag}, one brings the solution constructed in the proof of Theorem \ref{intermediaire} to rest at time $t=0$ and $t=T$ with the (Eulerian) local controllability of (\ref{KdV})-(\ref{KdV0}). One notices that while the regularity of the solution constructed in Theorem \ref{intermediaire} is sufficient to define the flow pointwisely, the usual $L^2$ local controllability of (\ref{KdV})-(\ref{KdV0}) is not. We therefore use the local controllability result of Zhang. To state the result, let us consider a non-zero initial data
\begin{equation}\label{KdVy0}
y(x,0)=y_0(x), \quad x\in [0,L]
\end{equation}
and let us denote the set of equations (\ref{KdV}), (\ref{KdVu}), (\ref{KdVv}), (\ref{KdVw}) and (\ref{KdVy0}) by (\ref{KdV})-(\ref{KdVy0}). In \cite{Zhang}, the following was proven. 

\begin{thm}\label{contreg}
Let $T>0$ and $s\geq 0$ be given and $[0,L]\subset (\alpha_1,\beta_1)$. Suppose that 
\[
w \equiv w(x,t) \in C^\infty((\alpha_1,\beta_1)\times (-\epsilon,T+\epsilon)),
\]
for some $\epsilon>0$, satisfies 
\[
w_t+w_x+w_{xxx}+ww_x=0, \qquad (x,t)\in (\alpha_1,\beta_1)\times (-\epsilon,T+\epsilon).
\]
Then there exists $\delta>0$ such that for any $y_0, y_T \in H^s(0,L)$ satisfying 
\[
\|y_0-w(.,0)\|_{H^s(0,L)}\leq \delta \quad and \quad \| y_T-w(.,T)\|_{H^s(0,L)} \leq \delta,  
\]
one can find control inputs $(u,v,w)\in H^{\frac{s+1}{3}}(0,T) \times H^{\frac{s+1}{3}}(0,T) \times H^{\frac{s}{3}}(0,T)$ such that (\ref{KdV})-(\ref{KdVy0}) has a solution 
\[
y\in C([0,T];H^s(0,L))\cap L^2((0,T);H^{s+1}(0,L)),
\]
satisfying 
\[
y(x,0)=y_0(x) \quad \textrm{and } \quad y(x,T)=y_T(x),
\]
on the interval $(0,L)$.
\end{thm}

Using Theorem \ref{contreg} with $s=2$ to connect the solution of Theorem \ref{intermediaire} to rest is sufficient to define the flow pointwisely, hence the statement of Theorem \ref{Lag} in $C([0,T];H^2(0,L))$. 

During the control phases, from the lack of maximum principle for the KdV equation, one cannot insure that the flow does not exit the domain $[0,L]$ from the left. Therefore, we use a stability estimate of the controls of Theorem \ref{contreg} with respect to the initial data to estimate the flow during the control phases. The following corollary is a consequence of the results in \cite{Zhang}.

\begin{cor}\label{stabilitecontrole}
Let $T>0$ and $s\geq 0$. Then, there exists $\delta>0$ such that for any $y_0, y_T\in H^s(0,1)$ satisfying 
\[
\|y_0\| < \delta \quad \textrm{and} \quad \|y_T\| < \delta
\]
there exists control inputs $(u,v,w)\in H^{\frac{s+1}{3}}(0,T) \times H^{\frac{s+1}{3}}(0,T) \times H^{\frac{s}{3}}(0,T)$ such that (\ref{KdV})-(\ref{KdVy0}) has a solution 
\[
y\in C([0,T];H^s(0,L))\cap L^2((0,T);H^{s+1}(0,L)),
\]
satisfying 
\[
y(x,0)=y_0(x) \quad \textrm{and } \quad y(x,T)=y_T(x),
\]
on the interval $(0,L)$. Moreover, 
\[
\|u\|_{H^{\frac{s+1}{3}}(0,T)}^2+\|v\|_{H^{\frac{s+1}{3}}(0,T)}^2+\|w\|_{H^{\frac{s}{3}}(0,T)}^2\leq c\left(\|y_0\|^2_{H^s(0,L)}+\|y_T\|^2_{H^s(0,L)}\right)
\]
where $c>0$ is independant of $y_0$ and $y_T$.
\end{cor}

Corollary \ref{stabilitecontrole} is the stability estimate of Theorem \ref{contreg} in the simpler case $w\equiv 0$ for which the Banach Fixed Point Theorem can be used. We comment the proof of the estimate in Section 2 for sake of completeness. 

Aside of Theorem \ref{contreg}, several other results of Eulerian controllability for (\ref{KdV})-(\ref{KdVy0}) are found in literature. When only the control $w$ is used ($u,v\equiv 0$), the local controllability around the equilibrium of (\ref{KdV})-(\ref{KdVy0}) was obtained by Rosier (\cite{Rosier}) if $L\notin \N$, that is, when the linearized equation around the equilibrium is exactly controllable. If $L\in \N$, there exists an unreachable state subspace for the linearized equation around the equilibrium. Using the power expansion method (\cite{CorBook}) to reach this subspace, the local controllability around the equilibrium of (\ref{KdV})-(\ref{KdVy0}) was obtained by Coron and Cr\'epeau (\cite{CoronCrepeau}), Cerpa (\cite{Cerpa}) and Cerpa and Cr\'epeau (\cite{CerpaCrepeau}) when the dimension of the unreachable states subspace is of dimension 1, 2 and of arbitrarily dimension, respectively. If one only uses $u$ as a control ($v,w\equiv 0$), then (\ref{KdV})-(\ref{KdVy0}) is locally controllable to zero (\cite{GG}). If one only uses $v$ ($u,w\equiv 0$), then (\ref{KdV})-(\ref{KdVy0}) is locally controllable around the equilibrium if $L$ doesn't belong to a countable set of critical lengths $\O$ (\cite{Glass}). If one only uses $v$ and $w$ ($u\equiv 0$), then it was shown in \cite{Rosier} that the system is small-time controllable. Any other combination of two controls also leads to the small-time controllability (\cite{GG}). Good surveys on the small-time eulerian controllability and stability of (\ref{KdV})-(\ref{KdVy0}) can be found in \cite{RZ,CerpaSurvey}

It is important to note here that, exception made of Theorem \ref{contreg}, none of the previously mentioned results of small-time Eulerian controllability of (\ref{KdV})-(\ref{KdVy0}) allows to define (\ref{flux}) pointwisely, the regularity of the solution of (\ref{KdV})-(\ref{KdVy0}) obtained with these results being at most in $C([0,T];L^2(0,L))\cap L^2([0,T];H^1(0,L))$. 

One finds in the literature two results of global Eulerian controllability for the KdV equation. Rosier proved in \cite{RosierGlobal} that (\ref{KdV})-(\ref{KdVy0}) is globally controllable, that is, that there are no smallness restrictions on the initial or final data. However, the minimal time of controllability $T>0$ depends on the initial and final data and may be large. By considering 
\begin{equation}
y_t+y_x+y_{xxx}+yy_x=a(t), \, x\in [0,L], t\in [0,T], \label{KdVChap}
\end{equation}
instead of (\ref{KdV}), Chapouly proved in \cite{Chapouly} the small-time global Eulerian controllability of (\ref{KdVChap}), (\ref{KdVu}), (\ref{KdVv}), (\ref{KdVw}) and (\ref{KdVy0}) using the controllability of the non viscous Burgers equation (\cite{ChapoulyVisqueuse}) where $a(t)$ is used as a fourth control. 

To date, the following challenging open problem still holds.
\begin{open}
Let $T,L>0$. For any $y_0\in L^2(0,L)$ and $y_T\in L^2(0,L)$, does there exist $u(t),v(t)$ and $w(t)\in \RR$ such that the solution $y$ of (\ref{KdV})-(\ref{KdVy0}) satisfies 
\[
y(x,T)=y_T(x)?
\]
\end{open}

Few results on Lagrangian controllability are found in the literature. Glass and Horsin showed for the 2-D (\cite{2dfluid}) and the 3-D (\cite{3dfluid}) Euler equation that, for two given smooth contractible sets of particles surrounding the same volume of fluids and any initial velocity field, it is possible to find a boundary control and a time interval such that the corresponding solution of the Euler equation makes the first set reaches approximately the second. Horsin proved in \cite{HorsinH} for the heat equation posed on $[0,L], L>0$ that if one consider any two closed interval of $[0,L]$, then there exists a boundary control such that the flow induced by the solution of the heat equation starting from the first interval reaches the second. In the same article, he proved in the case of a radial domain (respectively a convex domain) of higher dimension, one can move two regular closed sets with the flow induced by minus the gradient of the solution by a control action on a part of the domain in 
arbitrarily small-time (respectively sufficiently large time). Finally, Horsin proved, in the case of the viscous Burgers equation that the previously mentioned result holds locally, that is, for two intervals not too far apart (\cite{HorsinB}). 

Considering the controllability of a PDE written in Lagrangian coordinates, one finds the work of Rosier on the Korteweg-de Vries equation written in Lagrangian coordinates with wave-maker controls. He proved in \cite{RosierLag} the local controllability around regular trajectories.

One notes that the global controllability results of nonlinear equations are usually obtained by considering either the linear or nonlinear part of the equation as a perturbation. Fabre proved the approximate controllability of variations of the Navier-Stokes equation using the latter approach by truncating the nonlinearity (\cite{Fabre}). Still using the latter method, Fernandez-Cara and Zuazua \cite{FerZua} proved that, for the semilinear heat equation, the equation is null controllable if the nonlinear term is of controlled growth (see also \cite{Unsolved}, Problem 5.5, for a review and open problems of exact controllability of the semi-linear wave equation). Finally, Lions and Zuazua proved the controllability for some fluid systems using a Galerkin's approximation (\cite{LionsZua}). 

When considering the former approach, one may mimic the following finite dimensional result. Consider the finite dimensional system $y'=F(y)+Bu$ where $F$ is quadratic ($F(\lambda y)=\lambda^2F(y)$) and assume that there exists a trajectory $(\overline{y},\overline{u})$, satisfying $\overline{y}(0)=\overline{y}(T)=0$, of the system such that the linearized system around this trajectory is controllable. Then, by performing a scaling, one can show that $y'=F(y)+Ay+Bu$ is globally controllable $\forall A\in \LL(\RR^m;\RR^n)$ (\cite{CoronNew}). It is by using this technique that Chapouly proved the that the small-time global Eulerian controllability of (\ref{KdVChap})-(\ref{KdVy0}) (\cite{Chapouly}). This technique was also used by Coron and Glass to show, respectively, the global Eulerian controllability of the 2-D and 3-D Euler equation (\cite{CoronEuler,GlassEuler}). Coron \cite{CoronNavier} and Coron and Fursikov \cite{Fursikov} proved the global Eulerian controllability of the 2-D Navier-Stokes equations, 
using the global Eulerian controllability of the 2-D Euler equation, in the case where the whole boundary is used to control the interior or in the case of a Navier slip boundary condition. 

One remarks that if the result holds in finite dimension for any linear operator $A\in \LL(\RR^m;\RR^n)$, the presence of high order derivatives in the linear term and boundary layers issues may prevent one to apply this result in the infinite dimension framework. The first example where the global Eulerian controllability was obtain despite the presence of a boundary layer is due to Marbach, using the Hopf-Cole transformation and the maximum principle to show the small-time global null controllability of viscous Burgers equation (\cite{Marbach}).

The novelty of this paper is that we make full use of both the linearity and the nonlinearity of the KdV equation to obtain Theorem \ref{Lag}, as solitons don't exists if the linear or nonlinear term is dropped from (\ref{KdVHirota}). It is thus, to our knowledge, the first global controllability result obtained for a nonlinear equation without considering the linear or nonlinear part as a perturbation. 

The outline of the paper is the following. In Section 2, we state the well-posedness results for the linear and nonlinear KdV equation and review Corollary \ref{stabilitecontrole}. Section 3 is devoted to the proof of the main result.

\section{Well-posedness and regular controls}

\subsection{Well-posedness of the KdV equation}  

Consider the linear KdV equation

\begin{equation}\label{linwell}
 \left\{
 \begin{array}{l l}
y_t+y_x+y_{xxx}=0 , & x\in [0,L], \quad t\in [0,T],  \\ 
y(0,t)=u(t), y(L,t)=v(t), & t\in (0,T),\\
y_x(L,t)=w(t), & t\in (0,T), \\
y(x,0)=y_0(x), & x\in (0,L). 
\end{array}
\right.
\end{equation}

The well-posedness of (\ref{linwell}) was obtained in \cite{BonaZhang}, exhibiting the smoothing effects of the solution $y$ with respect to the initial value and the boundary data.
\begin{thm}
Let $L,T>0$. Let $y_0\in H^s(0,L)$ and $(u,v,w) \in H^{\frac{s+1}{3}}(0,T) \times H^{\frac{s+1}{3}}(0,T) \times H^{\frac{s}{3}}(0,T)$. Then, the problem (\ref{linwell}) has a unique solution in 
\[
 C([0,T];H^s(0,L))\cap L^2((0,T);H^{s+1}(0,L)). 
\]
Moreover, there exists $C>0$ such that 
\begin{align*}
 \|y\|_{C(0,T;H^s(0,L))\cap L^2((0,T);H^{s+1}(0,L))}^2\leq & C\left(  \|y_0\|_{H^{s}(0,L)}^2+\|u\|_{H^{s/3}(0,T)}^2 \right. \\  &  \quad  \left.+\|v\|_{H^{s/3}(0,T)}^2+\|w\|_{H^{s/3}(0,T)}^2\right).
\end{align*}
\end{thm}

Let us state the global well-posedness of the nonlinear KdV equation (\ref{KdV})-(\ref{KdVy0}) obtained in \cite{BonaZhang} (we refer to \cite{Faminskii} for a sharper result on the compatibility conditions). In order to state the well-posedness result, one needs to consider the $s$-compatibility conditions.
\begin{defn}[$s$-compatibility conditions]
Let $T,L,s>0$. A four-tuple $(\tilde{y},u,v,w)\in H^s(0,L)\times H^{(s+1)/3}(0,T)\times H^{(s+1)/3}(0,T)\times H^{s/3}(0,T)$ is said to be $s$-compatible if 
\begin{equation}\label{vacLG}
\tilde{y}_k(0)=u^{(k)}(0), \quad \tilde{y}_k(L)=w^{(k)}(0),
\end{equation}
hold for: 
\begin{enumerate}
\item $k=0,...,\floor{s/3}-1$ \textrm{when }$s-3 \floor{s/3}\leq 1/2$;
\item $k=0,...,\floor{s/3}$ \textrm{when }$3/2\geq s-3 \floor{s/3} >1/2$;
\end{enumerate}
and 
\begin{equation}
\tilde{y}_k(0)=u^{(k)}(0), \quad \tilde{y}_k(L)=w^{(k)}(0),\tilde{y}'_k(L)=h_3^{(k)}(0),
\end{equation}
holds for $k=0,...,\floor{s/3}$ when $s-3 \floor{s/3} >3/2$, where 
\begin{equation*}
\begin{cases}
\tilde{y_0}(x):=\tilde{y}(x), &\\ 
\tilde{y_k}(x):=-\tilde{y}_{k-1}^{(3)}(x)-\tilde{y}_{k-1}'(x)- \displaystyle \sum_{j=0}^{k-1} (\tilde{y}_j(x)\tilde{y}_{k-j-1}(x))', & k\in \NN. 
\end{cases}
\end{equation*}
We assume that \eqref{vacLG} is vacuous if $\floor{s/3}-1<0$.
\end{defn}

Let $\epsilon>0$ and 
\begin{equation*}
\mu_1(s):=
\begin{cases}
\epsilon+(5s+9)/18 & \textrm{ if } 0\leq s\leq 3, \\
(s+1)/3 & \textrm{ if } 3\leq s, \\
\end{cases}
\end{equation*}
\begin{equation*}
\mu_2(s):=
\begin{cases}
\epsilon+(5s+3)/18 & \textrm{ if } 0\leq s\leq 3, \\
(s+1)/3 & \textrm{ if } 3\leq s. \\
\end{cases}
\end{equation*}
Well-posedness follows from \cite[Theorem 1.3, p.1396]{BonaZhang} (see \cite{BonaZhang2} for the well-posedness for any $s>-1$),
\begin{thm}
For any $s\geq 0$, for any $T,L>0$ and for any $s$-compatible $(\tilde{y},u,v,w)\in H^s(0,L)\times H^{\mu_1(s)}(0,T)\times H^{\mu_1(s)}(0,T)\times H^{\mu_2(s)}(0,T)$, (\ref{KdV})-(\ref{KdV0}) is well-posed in 
\begin{equation*}
C([0,T];H^s(0,L))\bigcap L^2([0,T];H^{s+1}(0,L))
\end{equation*}
\end{thm}

\subsection{Regular controls}

Consider the linear KdV equation defined on the real line
\begin{equation}\label{cnt1}
\left\{
\begin{array}{ll}
z_t+z_x+z_{xxx}=0, &(x,t)\in \RR^2 \\
z(x,0)=z_0(x), &x\in \RR.
\end{array}
\right.
\end{equation}
The initial value control problem was solved for (\ref{cnt1}) in \cite{Zhang}
\begin{thm}[\cite{Zhang}, Theorem 3.1, p.554]\label{lintraj}
Let $s\geq 0$ and $T>0$ be given. There exists a bounded linear operator $G: H^s(0,L) \times H^s(0,L) \rightarrow H^s(\RR)$ such that for any $y_0, y_T \in H^s(0,L)$ if one chooses $z_0=G(y_0,y_T)\in H^s(\RR)$, then the corresponding solution $z$ of (\ref{cnt1}) satisfies 
\[
z(x,0)=y_0(x), \qquad z(x,T)=y_T(x),
\]
on the interval $(0,L)$ and 
\begin{equation}\label{IBVPestimate}
\|z_0\|_{H^s(\RR)}\leq c\left(\|y_0\|_{H^s(0,L)}+\|y_T\|_{H^s(0,L)}\right)
\end{equation}
where $c>0$ is independant of $y_0$ and $y_T$. 
\end{thm}
In fact, it was proven that the solution $z$ constructed in Theorem \ref{lintraj} is $C^{\infty}(\RR\times (0,T))$. 

A direct corollary of Theorem \ref{lintraj} is the exact controllability of the linear KdV equation by using the trace of $z$ as the boundary controls
\begin{equation}\label{cnt2}
\left\{
\begin{array}{l l}
y_t+y_x+y_{xxx}=0,  & x\in [0,L], \quad t\in [0,T],  \\ 
y(0,t)=u(t), y(L,t)=v(t) & t\in [0,T],  \\ 
y_x(L,t)=w(t) & t\in [0,T],  \\ 
y(x,0)=y_0(x),& x\in [0,L],  \\ 
\end{array}
\right.
\end{equation}
The corollary that we state here is slightly different than Corollary 3.4 in \cite[p. 559]{Zhang}
\begin{cor}\label{contlinreg} 
Let $s\geq 0$ and $T>0$ be given. For any $y_0, y_T \in H^s(0,L)$, there exists $(u,v,w)\in H^{\frac{s+1}{3}}(0,T) \times H^{\frac{s+1}{3}}(0,T) \times H^{\frac{s}{3}}(0,T)$, depending linearly on $y_0,y_T$, such that (\ref{cnt2}) has a solution
\[
y\in C([0,T];H^s(0,L))\cap L^2((0,T);H^{s+1}(0,L))
\]
satisfying 
\[
y(x,0)=y_0(x), \qquad y(x,T)=y_T(x),
\]
in the interval $(0,L)$. Moreover, there exists $C>0$, independant of $y_0$ and $y_T$ such that
\begin{equation}\label{estimecont}
\|u\|_{H^{\frac{s+1}{3}}(0,T)}^2+\|v\|_{H^{\frac{s+1}{3}}(0,T)}^2+\|w\|_{H^{\frac{s}{3}}(0,T)}^2
\leq C \left(\|y_0\|^2_{H^s(0,L)}+\|y_T\|^2_{H^s(0,L)}\right)
\end{equation}
\end{cor}
The operator $G$ constructed in the proof of Theorem \ref{lintraj} relies on the extension of the initial data $y_0$ of (\ref{cnt2}), defined on $[0,L]$, to the initial data of (\ref{cnt1}), defined on $\RR$, and of compact support. Since there exists infinitely many such extensions, the uniqueness of the controls of (\ref{cnt2}) is not guaranteed. The linearity of the controls $(u,v,w)$ with respect to $y_0,y_T$ stated in Corollary \ref{contlinreg} is obtained either by always choosing the same extension of the initial datas, either by considering the controls $(u,v,w)$ of minimal $H^{\frac{s+1}{3}}(0,T) \times H^{\frac{s+1}{3}}(0,T) \times H^{\frac{s}{3}}(0,T)$-norm since the projection is a linear operator . Moreover, the stability estimate (\ref{estimecont}) is obtained, thanks to the fact that the initial data $z_0$ constructed in the proof of Theorem \ref{lintraj} is compactly supported, from a sharp Kato smoothing effect of (\ref{cnt1}) (see for instance \cite{Kenig}) and from (\ref{IBVPestimate}). The regularity of the controls was stated in \cite{RZ}. Hence, there exists a linear continuous mapping from the initial data to the controls in the appropriate spaces. 

The last result needed to prove Corollary \ref{stabilitecontrole} is the following, which is the generalization of \cite[Proposition 4.1]{Rosier}

\begin{prop}\label{Prop41}
Let $s\geq 0$. Let $y\in L^2((0,T);H^{s+1}(0,L))$. Then, \mbox{$yy_x\in L^1((0,T);H^s(0,L))$} and the map $y\in L^2((0,T);H^{s+1}(0,L))\mapsto yy_x\in L^1((0,T);H^s(0,L))$ is continuous. Moreover, there exists $K>0$ such that, if $(y,z)\in L^2((0,T);H^{s+1}(0,L))^2$, then
\begin{eqnarray}\label{map}
\|yy_x-zz_x\|_{L^1((0,T);H^s(0,L))}\leq& K\left( \|y\|_{L^2((0,T);H^{s+1}(0,L))}+\|z\|_{L^2((0,T);H^{s+1}(0,L))}\right) \nonumber \\ 
& \cdot \|y-z\|_{L^2((0,T);H^{s+1}(0,L))}. 
\end{eqnarray}
\end{prop}

\beginpf 

First, consider the case $n\in \NN \cup \{0\}$. Let $y,z\in L^2((0,T);H^{n+1}(0,L))$. Using the Sobolev embedding of $H^{n+1}(0,L)$ in $L^\infty(0,L)$, we have
\begin{multline*}
\| yy_x-zz_x \|_{L^1((0,T);H^n(0,L))}\leq \int_0^T \|(y-z)y_x\|_{H^n(0,L)}
+\| z(y_x-z_x)\|_{H^n(0,L)}\, \textrm{dt} \\
\leq C  \int_0^T \|y-z\|_{L^\infty(0,L)}\|y_x\|_{H^n(0,L)} +\|y-z\|_{H^{n+1}(0,L)}\|y\|_{H^n(0,L)}  \\
 + \|z\|_{L^\infty(0,L)}\|y_x-z_x\|_{H^n(0,L)} +\|z\|_{H^{n+1}(0,L)}\|y-z\|_{H^n(0,L)} \,\textrm{dt}\\
\leq K\|y-z\|_{L^2((0,T);H^{n+1}(0,L))}
 \cdot (\|y\|_{L^2((0,T);H^{n+1}(0,L))}+\|z\|_{L^2((0,T);H^{n+1}(0,L))}) 
\end{multline*}
The result for $s\geq 0$ follows by interpolation (\cite{Tartar}).

\endpf

We obtain Corollary \ref{stabilitecontrole} by proving that the nonlinear equation is locally controllable with the Banach Fixed Point Theorem, the continuity of the controls with respect of the initial data being preserved by the continuity of the linear operators considered in the argument. The Banach Fixed Point Theorem argument to obtain the local controllability of a nonlinear equation from the controllability of the linearized equation is classical (see \cite{CorBook}).

\beginpf

Let $y_0, y_T\in H^s(0,L)$ such that $\| y_0 \|_{H^s(0,L)}\leq r$ and $\| y_T \|_{H^s(0,L)}\leq r$ with $r>0$ to be chosen later on. Consider $y^1, y^2, y^3$ the solutions of the following problems
\begin{equation}\label{lin1}
\left\{
\begin{array}{ll}
y^1_t+y^1_x+y^1_{xxx}=0,  & x\in [0,L], \quad t\in [0,T],  \\ 
y^1(0,t)=0, & t\in [0,T],  \\ 
y^1(L,t)=0,  & t\in [0,T],  \\ 
y^1_x(L,t)=0,  & t\in [0,T],  \\ 
y^1(x,0)=y_0(x),  & x\in [0,L],  \\ 
\end{array}
\right.
\end{equation}
\begin{equation}\label{lin2}
\left\{
\begin{array}{ll}
y^2_t+y^2_x+y^2_{xxx}=0,  & x\in [0,L], \quad t\in [0,T],  \\ 
y^2(0,t)=u(t), & t\in [0,T],  \\ 
y^2(L,t)=v(t), & t\in [0,T],  \\ 
y^2_x(L,t)=w(t), & t\in [0,T],  \\ 
y^2(x,0)=0,& x\in [0,L],  \\ 
\end{array}
\right.
\end{equation}
\begin{equation*}
\left\{
\begin{array}{ll}
y^3_t+y^3_x+y^3_{xxx}=f, & x\in [0,L], \quad t\in [0,T],  \\ 
y^3(0,t)=0, & t\in [0,T],  \\ 
y^3(L,t)=0, & t\in [0,T],  \\ 
y^3_x(L,t)=0, & t\in [0,T],  \\ 
y^3(x,0)=0, & x\in [0,L].  \\ 
\end{array}
\right.
\end{equation*}

Consider the continuous maps, 
\begin{equation*}
\begin{array}{rccl}
 \psi_1  : & H^{\frac{s+1}{3}}(0,T) \times H^{\frac{s+1}{3}}(0,T) \times H^{\frac{s}{3}}(0,T) & \rightarrow &  L^2((0,T);H^{s+1}(0,L))\\ 
 & (u,v,w) & \mapsto& y^2, 
\end{array}
\end{equation*}
and 
\begin{equation*}
\begin{array}{rccl}
 \psi_2 : &L^1((0,T);H^{s}(0,L))& \rightarrow & L^2((0,T);H^{s+1}(0,L))\\ 
 & f & \mapsto& y3, 
\end{array}
\end{equation*}
Let
\begin{equation*}
\begin{array}{rccl}
 \Gamma : &H^{s}(0,L)& \rightarrow &H^{\frac{s+1}{3}}(0,T) \times H^{\frac{s+1}{3}}(0,T) \times H^{\frac{s}{3}}(0,T) \\ 
 &y_T& \mapsto& (u,v,w), 
\end{array}
\end{equation*}
be the continuous map associating to $y_T$, the controls $(u,v,w) \in H^{\frac{s+1}{3}}(0,T) \times H^{\frac{s+1}{3}}(0,T) \times H^{\frac{s}{3}}(0,T)$ given by Corollary \ref{contlinreg}, such that $y^2$, the solution of the backward equation (\ref{lin2}) starting from $y^2(.,T)=y_T$, reaches $y^2(.,0)=0$. Let
\begin{equation*}
\begin{array}{rccl}
 F : &L^2((0,T);H^{s+1}(0,L))& \rightarrow &L^2((0,T);H^{s+1}(0,L))\\
  &y& \mapsto& F(y),
  \end{array}
\end{equation*}
with 
\begin{equation*}
 F(y):=y^1+\psi_1 \circ \Gamma \left( y_T-y^1(.,T)+\psi_2(yy_x)(.,T)\right) + \psi_2(-yy_x).
\end{equation*}
The map $F$ was constructed such that it is well-defined, continuous and that every fixed point of $F$ is a solution of (\ref{KdV})-(\ref{KdVy0}) satisfying $y(.,T)=y_T$. Therefore, it is sufficient to prove that, for a closed ball $B(0,R)\subset L^2((0,T);H^{s+1}(0,L))$, $F(B(0,R))\subset B(0,R)$, and that there exists $C\in (0,1)$ such that, $\forall (y,z)\in B(0,R)^2$,
\[
 \|F(y)-F(z)\|_{L^2((0,T);H^{s+1}(0,L))}\leq C\|y-z\|_{L^2((0,T);H^{s+1}(0,L))},
\]
to show the existence of a fixed point of $F$ by the Banach Fixed Point Theorem.

Let $K_1,K_\Gamma$ be the norms of $\psi_1, \Gamma$ and let $K_2,K_2'$ be the norm of $\psi_2$ in \newline $L^2((0,T);H^{s+1}(0,L))$ and $C((0,T);H^s(0,L))$ respectively. Furthermore, let $K$, be the constant from Proposition \ref{Prop41}, $C_1$ denotes the constant in the stability estimate for solutions of (\ref{lin1}),
\[
\|y^1\|_{L^2(0,T;H^{s+1}(0,L))}+\|y^1\|_{C(0,T;H^s(0,L))} \leq C_1 \|y_0\|_{H^s(0,L)}.
\]
Then, for $y\in L^2((0,T);H^{s+1}(0,L))$, 
\begin{align*}
\|F(y)\|_{L^2((0,T);H^{s+1}(0,L))}\leq& \|y^1\|_{L^2(0,T;H^{s+1}(0,L))}+K_1K_\Gamma \left(\|y_T\|_{H^s(0,L)} \right. \\
& \left.+\|y^1(.,T)\|_{H^s(0,L)} +\|\psi_2(yy_x)(.,T)\|_{H^s(0,L)}\right)\\
&+\|\psi_2(-yy_x)\|_{L^2(0,T;H^{s+1}(0,L))} \\
\leq& C_1\|y_0\|_{H^s(0,L)}+K_1K_\Gamma \left(\|y_T\|_{H^s(0,L)} \right. \\
& \left. +\|y^1\|_{C(0,T);H^s(0,L))} +\|\psi_2(yy_x)\|_{C(0,T;H^s(0,L))}\right) \\
&+K_2 \|yy_x\|_{L^1(0,T;H^s(0,L))} \\
\leq& C_1\|y_0\|_{H^s(0,L)}+K_1K_\Gamma \left(\|y_T\|_{H^s(0,L)}  \right. \\
& \left. +C_1\|y_0\|_{H^s(0,L)}+K_2'\|yy_x\|_{L^1(0,T;H^s(0,L))}\right)\\ 
&+K_2\|yy_x\|_{L^1(0,T;H^s(0,L))} 
\end{align*}
\begin{align*}
\qquad \qquad \qquad \qquad \, \, \, \, \, \, \, \, \leq& C_1\|y_0\|_{H^s(0,L)}+K_1K_\Gamma \big(\|y_T\|_{H^s(0,L)} \\
&  +C_1\|y_0\|_{H^s(0,L)} +K_2'K\|y\|^2_{L^2(0,T;H^{s+1}(0,L))}\big)\\
&+K_2 K\|y\|^2_{L^2(0,T;H^{s+1}(0,L))} \\
\leq& C_1r+K_1K_\Gamma \left(r+C_1r+K_2'KR^2\right)+K_2KR^2,
\end{align*}
where all the constants are independent of $r$ and $R$. We impose on $r$ and $R$ that 
\begin{equation}\label{estm1}
C_1r+K_1K_\Gamma \left(r+C_1r+K_2'KR^2\right)+K_2KR^2<R.
\end{equation}
Furthermore, 
\begin{align*}
 \|F(y)-F(z)\|_{L^2((0,T);H^{s+1}(0,L))}=&\| \psi_1 \circ \Gamma \left( \psi_2(yy_x-zz_x)(.,T)\right)  \\
 &+ \psi_2(-yy_x+zz_x) \|_{L^2((0,T);H^{s+1}(0,L))}\\
 \leq& K_1K_\Gamma \| \psi_2(yy_x-zz_x)(.,T) \|_{H^s(0,L)} \\
 & +  \| \psi_2(-yy_x+zz_x) \|_{L^2((0,T);H^{s+1}(0,L))}\\
 \leq& K_1K_\Gamma \| \psi_2(yy_x-zz_x) \|_{C(0,T;H^s(0,L))} \\ 
 &+ K_2 \| -yy_x+zz_x \|_{L^1((0,T);H^s(0,L))}\\
 \leq& K_1K_\Gamma K_2'  \| yy_x-zz_x \|_{L^1((0,T);H^s(0,L))} \\
 &+ K_2 \| -yy_x+zz_x \|_{L^1((0,T);H^2(0,L))}\\
 \leq& 2R(K_1K_\Gamma K_2'K+K_2K) \\
 &\cdot \|y-z\|_{L^2((0,T);H^{s+1}(0,L))}.
\end{align*}
We obtain that $F$ is a contraction by choosing $R>0$ small enough so that
\begin{equation*}
2R(K_1K_\Gamma K_2'K+K_2K)<1.
\end{equation*}
By taking into account (\ref{estm1}), we then choose
\[
r=\dfrac{R}{2(C_1+K_1K_\Gamma+K_1K_\Gamma C_1)}.
\]
Thus there exists a fixed point of $F$. 

\endpf

\section{Lagrangian controllability}

First, let us prove Theorem \ref{Lag} assuming Theorem \ref{intermediaire}.
\newline

\beginpf 

Let $L,T>0$, $(\epsilon_1,\epsilon_2)\in (0,T/2)^2$ and $\gamma>0$. For every $(a,b)\in \RR^2$, such that $a<0<L<b$ and $b-a>L+\gamma$, let us denote by $\widetilde{y_0}$ an extension of $y_0$ in $H^2(a,b)$ with homogeneous Dirichlet boundary conditions. It is well-known that this extension can be chosen so that there exists $C(a,b)>0$ such that 
\begin{equation}
\|\widetilde{y_0}\|_{H^2(a,b)}\leq C(a,b) \|y_0\|_{H^2(0,L)}.
\end{equation}
By the classical stability estimate of solutions of (\ref{KdV})-(\ref{KdV0}), by Corollary \ref{stabilitecontrole} and the Sobolev embedding of $H^2(a,b)$ in $L^\infty(a,b)$, there exists a solution $y$ solution of (\ref{KdV})-(\ref{KdV0}) starting from 0 to $y_0$ such that there exists $C_2>0$ such that 
\begin{equation}
\|y\|_{L^{\infty}((a,b)\times (0,T))} \leq C_2 \|y_0\|_{H^2(0,L)}.
\end{equation}
Therefore, let $y^2$ denotes the solution given by Theorem \ref{intermediaire} such that $\|y^2(.,\epsilon_1)\|_{H^2(0,L)}$ and $\|y^2(.,\epsilon_2)\|_{H^2(0,L)}$ are smaller than the tolerance $\delta$ stated in Corollary \ref{stabilitecontrole}, \newline $C_2 \|y^2(.,\epsilon_i)\|_{H_2(0,L)} \leq \gamma/2$, for $i=1,2$, and such that 
\[
\Phi(x,T)\geq L+\delta, \forall x\in [0,L].
\] 
Thus, Corollary \ref{stabilitecontrole} implies that there exists $(u^1,v^1,w^1)\in H^1(0,T)\times H^1(0,T)\times H^{\frac{2}{3}}(0,T) $ such that the solution $y^1\in C([0,\epsilon_1];H^2(a,b))$ of
\begin{equation*}
 \left\{
 \begin{array}{l l}
y^1_t+y^1_x+y^1_{xxx}+y^1y^1_x=0 , & x\in [a,b] \quad t\in [0,\epsilon_1],  \\ 
y^1(a,t)=u^1(t) , & t\in [0,\epsilon_1],  \\
y^1(b,t)=v^1(t) , & t\in [0,\epsilon_1],  \\
y^1_x(b,t)=w^1(t) , & t\in [0,\epsilon_1], \\
y^1(x,0)=0, & x\in [a,b], 
\end{array}
\right.
\end{equation*}
satisfies $y^1(.,\epsilon_1)=\widetilde{y^2}(.,\epsilon_1)$ and that there exists $(u^3,v^3,w^3)\in H^1(0,T)\times H^1(0,T)\times H^{\frac{2}{3}}(0,T) $ such that $y^3\in C([T-\epsilon_2,T];H^2(a,b))$, solution of 
\begin{equation*}
 \left\{
 \begin{array}{l l}
y^3_t+y^3_x+y^3_{xxx}+y^3y^3_x=0 , & x\in [a,b] \quad t\in [T-\epsilon_2,T],  \\ 
y^3(a,t)=u^3(t) , & t\in [T-\epsilon_2,T],  \\
y^3(b,t)=v^3(t) , & t\in [T-\epsilon_2,T],  \\
y^3_x(b,t)=w^3(t) , & t\in [T-\epsilon_2,T], \\
y^3(x,T-\epsilon_2)=\widetilde{y^2}(.,\epsilon_2), & x\in [a,b], 
\end{array}
\right.
\end{equation*}
satisfies $y^3(.,T)=0$. 

Let 
\begin{equation*}
  u(t):=\left\{
 \begin{array}{l l}
u^1(t), & t\in(0,\epsilon_1), \\
y^2(0,t), & t\in(\epsilon_1,T-\epsilon_2), \\
u^3(t), & t\in(T-\epsilon_2,T), 
\end{array}
\right.
\end{equation*}
\begin{equation*}
  v(t):=\left\{
 \begin{array}{l l}
v^1(t), & t\in(0,\epsilon_1), \\
y^2(L,t), & t\in(\epsilon_1,T-\epsilon_2), \\
v^3(t), & t\in(T-\epsilon_2,T), 
\end{array}
\right.
\end{equation*}
\begin{equation*}
  w(t):=\left\{
 \begin{array}{l l}
w^1(t), & t\in(0,\epsilon_1), \\
y^2_x(L,t), & t\in(\epsilon_1,T-\epsilon_2), \\
w^3(t), & t\in(T-\epsilon_2,T). 
\end{array}
\right.
\end{equation*}
Then, $y$ solution of
\begin{equation*}
 \left\{
 \begin{array}{l l}
y_t+y_x+y_{xxx}+yy_x=0 , & x\in [0,L] \quad t\in [0,T],  \\ 
y(0,t)=u(t) , & t\in [0,T],  \\
y(L,t)=v(t) , & t\in [0,T],  \\
y_x(L,t)=w(t) , & t\in [0,T], \\
y(x,0)=0, & x\in [0,L], 
\end{array}
\right.
\end{equation*}
belongs to $C([0,T];H^2(0,L))$ and, by construction of $y$, we obtain Theorem \ref{Lag}.

\endpf

We now conclude with the proof of Theorem \ref{Lag}
\newline

\beginpf

Let $T,L,\delta>0$ and $(\epsilon_1,\epsilon_2)\in (0,T/2)^2$. Let $\alpha_1>0$ and $\epsilon>0$. We define 
\[N:=\ceil[\bigg]{4L \alpha_1^2 / \ln\left(\sqrt{2\alpha_1}\left(1+\sqrt{1-\frac{1}{2\alpha_1}}\right)\right)}\]
where $\ceil{x}$ is the ceiling function. Let
\begin{equation}
0<\alpha_N <...< \alpha_1. \label{a}
\end{equation}
where 
\begin{equation}\label{diff}
 \alpha_1-\alpha_N=\epsilon,
\end{equation}
 and, for $i=1,...,N$,
\begin{align}
 s_i:=&L-(\alpha_1-\epsilon)^2(T-\epsilon_2)+\frac{N-i+1}{N+1}\left(-\alpha_1^2\epsilon_1-L+(\alpha_1-\epsilon)^2(T-\epsilon_2)\right) \nonumber \\
 =&-\alpha_1^2\epsilon_1\frac{N+1-i}{N+1}+\frac{i}{N+1}\left(L-(\alpha_1-\epsilon)^2(T-\epsilon_2)\right)\label{s}
\end{align}
Finally, let $G:=\left( \ln F\right)_{xx}$ and   
\[
F=1+\sum_{n=1}^N\sum_{C_{n}^{N}}a(i_1,...,i_n) \prod^{n}_{j=1} f_{i_j},
\]
where
\[
f_i(x,t)=\exp(-\alpha_i(x-s_i)+\alpha_i^3t).
\]
The function $G$ corresponds to $N$ solitons, ordered, for $t\in [0,T]$, from the right with the fastest soliton, associated to $\alpha_1$, to left with the slowest, associated to $\alpha_N$. Moreover, they were constructed so that, for $\alpha_1$ large, the $N$ solitons are located to the left of the interval $[0,L]$ for time $t\in [0,\epsilon_1]$, that they pass inside the domain during the time interval $(\epsilon_1,T-\epsilon_2)$ and are located to the right of $[0,L]$ for $t\in (T-\epsilon_2,T)$. We now prove that, for $\alpha_1$ sufficiently large, the function $G$ fulfills the requirements of Theorem \ref{intermediaire}. 

One notes from (\ref{s}) that, if 
\begin{equation}\label{speed1}
 (\alpha_1-\epsilon)^2>L/(T-\epsilon_2), 
\end{equation}
a natural condition on the speed of the slowest soliton for it to cross the domain $[0,L]$ in the required time, then,
\[
 s_N < \ldots < s_1 < 0.
\]
Assumption (\ref{speed1}) on $\alpha_1$ is therefore made for the rest of the proof. 

Let us show that there exists $\alpha_1$ sufficiently large such that $\|G(.,t)\|^2_{H^2} \leq \delta, t\in [0,\epsilon_1]$. 
The expression of $G$ and its first two derivatives takes the form,
\begin{align}
G&=\dfrac{FF_{xx}-(F_x)^2}{F^2}, \label{es1} \\
G_x&=\dfrac{F^2F_{xxx}-3FF_xF_{xx}+2(F_x)^3}{F^3}, \label{es2} \\
G_{xx}&=\dfrac{-4F^2F_xF_{xxx}+F^3F_{xxxx}+12F(F_x)^2F_{xx}-3(FF_{xx})^2-6(F_x)^4}{F^4}. \label{es3}
\end{align}
Thus, we have the bounds, thanks to the fact that $F(x,t) \geq 1, \, \forall (x,t)\in \RR^2$, 
\begin{align*}
|G|&\leq |F||F_{xx}|+|F_x|^2, \\
|G_x|&\leq |F|^2|F_{xxx}|+3|F||F_x||F_{xx}|+2|F_x|^3, \\
|G_{xx}|&\leq 4|F|^2|F_x||F_{xxx}|+|F|^3|F_{xxxx}|+12|F||F_x|^2|F_{xx}|+3|F|^2|F_{xx}|^2+6|F_x|^4.
\end{align*}
By noting that the derivatives of $F$ are given, for $k\in \NN$, by
\[
\dfrac{d^kF}{dx^k}=\sum_{n=1}^N\sum_{C_{n}^{N}}(-1)^k(\alpha_{i_1}+...+\alpha_{i_n})^k a(i_1,...,i_n) \prod^{n}_{j=1} f_{i_j}, 
\]
we see, taking into account (\ref{a}), (\ref{s}), that for $t\in [0,\epsilon_1]$, 
\[
\left\| \dfrac{d^kF}{dx^k}(x,t)\right\|_{L^\infty(0,L)}= \left|\dfrac{d^kF}{dx^k}(0,\epsilon_1)\right|.
\]
Therefore, taking into account that $a(i_1,...,i_n)$ are at least bounded by (\ref{diff}), if  
\begin{equation} \label{cond1}
 s_i+\alpha_i^2\epsilon_1 <0,
\end{equation}
hold, that is a condition insuring that all the solitons are located to the left of $[0,L]$, then there exists $\alpha_1$ sufficiently large such that $\|G(.,t)\|^2_{H^2(0,L)} \leq \delta, t\in [0,\epsilon_1]$. Or, 
\begin{align*}
\alpha_i^2\epsilon_1+s_i <& \alpha_1^2\epsilon_1+s_i \\
 =& \frac{i}{N+1}(L-(\alpha_1-\epsilon)^2(T-\epsilon_2)+\alpha_1^2\epsilon_1) \\
  =& \frac{i}{N+1}(L-\alpha_1^2(T-\epsilon_1-\epsilon_2)+2\alpha_1\epsilon(T-\epsilon_2)-\epsilon^2(T-\epsilon_2)) \\
 <&0,
\end{align*}
the last line holding for sufficiently large $\alpha_1$. Thus, (\ref{cond1}), a condition insuring that the $N$ solitons are located in $x<0$ for $t\in [0,\epsilon_1]$, holds.

We now prove that for $\alpha_1$ sufficiently large, $\|G(.,t)\|^2_{H^2} \leq \delta, t\in [0,\epsilon_1]$. From (\ref{es1})-(\ref{es3}), one sees that, when the $f_i$'s are large, the leading term $(f_1\cdots f_N)^{j+2}$ of $G^{(j)}$, for $j=0,1,2$, is only present to the denominator. Therefore, one has that if 
\begin{equation}\label{cond2}
 s_i-L+\alpha_i^2(T-\epsilon_2)>0,
\end{equation}
then there exists $\alpha_1$ such that $\|G(.,t)\|^2_{H^2(0,L)} \leq \delta, t\in [T-\epsilon_2,T]$. Or, 
\begin{align*}
  s_i-L+\alpha_i^2(T-\epsilon_2) >& s_i-L+\alpha_N^2(T-\epsilon_2) \\
  =& -\alpha_1^2\epsilon_1\frac{N+1-i}{N+1}+\frac{i}{N+1}\left(L-(\alpha_1-\epsilon)^2(T-\epsilon_2)\right) \\
  & -L + (\alpha_1-\epsilon)^2(T-\epsilon_2)\\
  =& -\alpha_1^2\epsilon_1 + (\alpha_1-\epsilon)^2(T-\epsilon_2) - L \\ 
  & +\frac{i}{N+1}(\alpha_1^2\epsilon_1+L-(\alpha_1-\epsilon)^2(T-\epsilon_2))\\
  =& \alpha_1^2(T-\epsilon_1-\epsilon_2)-2\epsilon\alpha_1(T-\epsilon_2)+\epsilon^2-L +\frac{i}{N+1}\\
  & \cdot (-\alpha_1^2(T-\epsilon_1-\epsilon_2)+2\epsilon\alpha_1(T-\epsilon_2)-\epsilon^2(T-\epsilon_2)+L) \\ 
  >&0,
\end{align*}
holds for sufficiently large $\alpha_1$. One remarks that (\ref{cond2}) represents the fact that the $N$ solitons are located in $x>L$ for $t\in [T-\epsilon_2,T]$. Moreover, combining (\ref{cond1}) and (\ref{cond2}) asks for $\alpha_i^2(T-\epsilon_1-\epsilon_2)>L$, a natural condition on the travelling speed of the $N$ solitons. Moreover, we have 
\[
L-\alpha_N^2(T-\epsilon_2)<s_N < \ldots < s_1 < -\alpha_1^2 \epsilon_1
\]

In order to prove that the flow $\Phi$ associated to $G$ satisfies $\Phi(x,T)\geq L, (x,t)\in (0,L)\times (T-\epsilon_2,T)$, we use a rough lower bound on each solitons of $G$ : a characteristic function of the same amplitude and the same width. Prior this estimate, let us first obtain rigorously that $G$ is greater than $N$ solitons. It is important to note here that, from the definition of $G$, we have that every coefficient in front of $f_{i_1}\dotsm f_{i_m}$, $1\leq m \leq 2N$ with, possibly, twice repeated indexes $i_j$, are positive. 

For a fixed time $t$ in $[0,T]$, in the neighbourhood of the first soliton $x\in ((3s_1+s_2)/4 +\alpha_1^2 t, -\alpha_1^2 \epsilon_1+\alpha_1^2 t)$, we have,

\begin{align*}
G(x,t) \geq& \dfrac{\alpha_1^2 f_1(x,t)}{\Big(1+\displaystyle{\sum_{n=1}^N\sum_{C_{n}^{N}}}a(i_1,...,i_n) \prod^{n}_{j=1} f_{i_j}(x,t)\Big)^2 } \\ 
 =& \dfrac{\alpha_1^2 f_1(x,t)}{\Big(1+\displaystyle{\sum_{n=1}^N\sum_{\substack{C_{n}^{N} \\ i_1 \neq 1}}}a(i_1,...,i_n) \prod^{n}_{j=1} f_{i_j}(x,t) + f_1(x,t)\Big)^2 } \\ 
\geq & \dfrac{\alpha_1^2 f_1(x,t)}{\Big(1+\displaystyle{\sum_{n=1}^N\sum_{\substack{C_{n}^{N} \\ i_1 \neq 1}}}a(i_1,...,i_n) \prod^{n}_{j=1} f_{i_j}\Big(\frac{3s_1+s_2}{4}+\alpha_1^2 t,t \Big) + f_1(x,t)\Big)^2 } 
\end{align*}
\begin{align*}
\qquad \qquad \quad \, \,=& \dfrac{1}{\Big(1+\displaystyle{\sum_{n=1}^N\sum_{\substack{C_{n}^{N} \\ i_1 \neq 1}}}a(i_1,...,i_n) \prod^{n}_{j=1} f_{i_j}\Big(\frac{3s_1+s_2}{4}+\alpha_1^2 t,t \Big)  \Big)^2}   \\
&  \cdot  \dfrac{\alpha_1^2 f_1(x,t)}{\Big(1+\Big( \frac{1}{1+\sum_{n=1}^N\sum_{\substack{C_{n}^{N} \\ i_1 \neq 1}}a(i_1,...,i_n) \prod^{n}_{j=1} f_{i_j}((3s_1+s_2)/4 +\alpha_1^2 t,t)  }\Big)f_1(x,t)\Big)^2 } \\
= & \dfrac{1}{\Big(1+\displaystyle{\sum_{n=1}^N\sum_{\substack{C_{n}^{N} \\ i_1 \neq 1}}}a(i_1,...,i_n) \prod^{n}_{j=1} f_{i_j}\Big(\frac{3s_1+s_2}{4}+\alpha_1^2 t,t \Big)  \Big)}   \\
& \cdot  \dfrac{\alpha_1^2\Bigg( \frac{1}{1+\sum_{n=1}^N\sum_{\substack{C_{n}^{N} \\ i_1 \neq 1}}a(i_1,...,i_n) \prod^{n}_{j=1} f_{i_j}((3s_1+s_2)/4 +\alpha_1^2 t,t)  }\Bigg)  f_1(x,t)}{\Big(1+\Big( \frac{1}{1+\sum_{n=1}^N\sum_{\substack{C_{n}^{N} \\ i_1 \neq 1}}a(i_1,...,i_n) \prod^{n}_{j=1} f_{i_j}((3s_1+s_2)/4 +\alpha_1^2 t,t)  }\Big) f_1(x,t)\Big)^2 }.
\end{align*}
Let
\[
A_1(t):=\dfrac{1}{\Big(1+\displaystyle{\sum_{n=1}^N\sum_{\substack{C_{n}^{N} \\ i_1 \neq 1}}}a(i_1,...,i_n) \prod^{n}_{j=1} f_{i_j}\Big(\frac{3s_1+s_2}{4}+\alpha_1^2 t,t\Big)  \Big)}. 
\]
Then, we obtain, in the region $x\in ((3s_1+s_2)/4 +\alpha_1^2 t, -\alpha_1^2 \epsilon_1+\alpha_1^2 t)$, 
\begin{equation}\label{soln1}
G(x,t)\geq \dfrac{\alpha_1^2A_1(t)}{4}\sech^2\left(\dfrac{-\alpha_1(x-\sigma_1(t)-\alpha_1^2 t)}{2}\right),
\end{equation}
a soliton of amplitude $\alpha_1^2A_1(t)/4$ of phase $\sigma_1(t):=s_1+\frac{1}{\alpha_1}\ln(A_1(t))$. 

For fixed $t$ in $[0,T]$, we consider, for the $k$-th soliton, $2 \leq k \leq N$, the neighbourhood 
\begin{equation*}
  (\xi^{-}_{k}(t),\xi^{+}_{k}(t))= \begin{cases}
  \left( \dfrac{2s_k+s_{k+1}}{3}+\alpha_k^2 t, \dfrac{2s_k+s_{k-1}}{3}+\alpha_k^2 t\right), \quad 2\leq k\leq N-1, \\
  \left( \dfrac{4s_N-s_{N-1}}{3}+\alpha_N^2 t, \dfrac{2s_N+s_{N-1}}{3}+\alpha_N^2 t\right), \quad  k=N.
 \end{cases}
\end{equation*}
One notices, from the definition of $s_N$, that $L-\alpha_N^2(T-\epsilon_2)<(4s_N-s_{N-1})/3<s_N$. 

Let 
 \[
 A_k(t)=\left(\dfrac{1}{1+\frac{1+\sum_{n=1}^N\sum_{_{N}C_{n,n\neq k-1,k}}a(i_1,...,i_n) \prod^{n}_{\substack{j=1 \\ i_j \leq k-1}} f_{i_n}(\xi^{+}_k(t),t)\prod^{n}_{\substack{j=1 \\ i_j > k-1}} f_{i_n}(\xi^{-}_k(t),t)}{a(1,...,k-1)f_1\dotsm f_{k-1}(\xi^{+}_k(t),t)} }\right),
\]
where $_{N}C_{n,n\neq k-1,k}$ denotes the usual $C_{n}^{N}$ but excluding the case where $i_1,...,i_n$ equal $1,...,k-1$ or $1,...,k$. Then, the same steps as in the case of the first soliton yield 
\begin{align*}
\hspace{-3cm}G(x,t) \geq& \dfrac{\alpha_k^2 a(1,...,k-1)a(1,...,k)f_1^2\dotsm f_{k-1}^2f_k(x,t)}{\left(1+\displaystyle{\sum_{n=1}^N\sum_{C_{n}^{N}}}a(i_1,...,i_n) \prod^{n}_{j=1} f_{i_j}(x,t)\right)^2 } 
\end{align*}
\begin{align*}
\qquad \qquad \, =& \dfrac{\alpha_k^2 \frac{a(1,...,k)}{a(1,...,k-1)}f_k(x,t)}{\left(\dfrac{1+\sum_{n=1}^N\sum_{\substack{ C_{n}^{N} \\ n\neq k-1,k}}a(i_1,...,i_n) \prod^{n}_{j=1} f_{i_j}(x,t)}{a(1,...,k-1)f_1\dotsm f_{k-1}(x,t)}+1+\frac{a(1,...,k)}{a(1,...,k-1)}f_k(x,t) \right)^2 } \\ 
\geq& \dfrac{\alpha_k^2 \frac{a(1,...,k)}{a(1,...,k-1)}f_k}{\left(A_k(t)^{-1}+\frac{a(1,...,k)}{a(1,...,k-1)}f_k \right)^2 } \\ 
 =& A_k(t)\dfrac{\alpha_k^2 A_k(t) \frac{a(1,...,k)}{a(1,...,k-1)}f_k}{\left(1+A_k(t)\frac{a(1,...,k)}{a(1,...,k-1)}f_k\right)^2},
 \end{align*}
Thus, for $x\in (\xi^{-}_{k}(t),\xi^{+}_{k}(t))$, we have
\begin{equation}\label{solnk}
G(x,t)\geq \dfrac{\alpha_k^2A_k(t)}{4}\sech^2\left(\dfrac{-\alpha_k(x-\sigma_k(t)-\alpha_k^2 t)}{2}\right),
\end{equation}
a soliton of amplitude $\alpha_k^2A_k(t)/4$ and of phase $\sigma_k:=s_k+\frac{1}{\alpha_k}\ln(\frac{A_k(t)a(1,...,k)}{a(1,...,k-1)})$. 

We prove the following for $A_k(t)$
\begin{lem}\label{lemAk}
 For $1\leq k \leq N$, we have, for $t\in [0,T]$, $A_k(t)\rightarrow 1$ as $\alpha_1 \rightarrow \infty$. 
\end{lem}

\beginpf

We first consider $k=1$. Let us remark that the term $f_1$ alone was removed from $A_1(t)$. Thus, let us show that $f_1f_2(\frac{2s_1+s_2}{3}+\alpha_1^2 t,t)$ converges to zero as $\alpha_1$ tends to infinity, since it is the biggest term in the sum.
\begin{align*}
 f_1f_2\left(\frac{2s_1+s_2}{3}+\alpha_1^2 t,t\right)=& \exp\left(-\alpha_1\left(\frac{2s_1+s_2}{3}+\alpha_1^2t-s_1-\alpha_1^2 t\right)\right) \\
 & \cdot   \exp\left(-\alpha_2\left(\frac{2s_1+s_2}{3}+\alpha_1^2t-s_2-\alpha_2^2 t\right)\right) \\
=& \exp\left((2\alpha_2-\alpha_1)\left(\frac{s_2-s_1}{3}\right)-\alpha_2(\alpha_1^2-\alpha_2^2)t \right). 
\end{align*}
Thus, since $s_2-s_1<0$ and, for large $\alpha_1$, $\alpha_1^2-\alpha_2^2>0$ and $2\alpha_2-\alpha_1>0$, since $\alpha_1-\alpha_N=\epsilon$, then the last expression converges to zero as $\alpha_1$ tends to infinity, and so do all the other terms in the sum of $A_1(t)$. 

We consider next the case $2\leq k\leq N-1$. We stress here that, since $f_1\dotsm f_{k-1}$ and $f_1\dotsm f_{k}$ were removed from the sum in $A_k(t)$, no terms of the form $f_{i_1}\dotsm f_{i_n}$, where $(i_1,\ldots,i_n) \in \{1,...,k-1\}^n$, appears at the numerator of 
\[
\frac{1+\sum_{n=1}^N\sum_{_{N}C_{n,n\neq k-1,k}}a(i_1,...,i_n) \prod^{n}_{\substack{j=1 \\ i_j \leq k-1}} f_{i_n}(\xi^{+}_k(t),t)\prod^{n}_{\substack{j=1 \\ i_j > k-1}} f_{i_n}(\xi^{-}_k(t),t)}{a(1,...,k-1)f_1\dotsm f_{k-1}(\xi^{+}_k(t),t)}. 
\]
Therefore, it is sufficient to show that the terms of the form $f_k(\xi_{k}^{-}(t),t)f_{k-1}^{-1}(\xi_k^{+}(t),t)$ and $f_k(\xi_{k}^{-}(t),t)f_{k+1}(\xi_k^{-}(t),t)$ tend to zero, as all the other terms will converge to zero as well. We have,
\begin{align*}
f_k(\xi_{k}^{-}(t),t)f_{k-1}^{-1}(\xi_k^{+}(t),t)=&\exp\left(-\alpha_k\left(\frac{-s_k+s_{k+1}}{3}\right) \right. \\ 
& \left.+\alpha_{k-1}\left(\frac{2s_k-2s_{k-1}}{3}\right)+(\alpha_k^2-\alpha_{k-1}^2)t\right) \\
=&\exp\left((\alpha_{k-1}-\alpha_k)\left(\frac{-s_k+s_{k+1}}{3}\right) \right. \\ 
& \left.+\alpha_{k-1}\left(\frac{3s_k-2s_{k-1}-s_{k+1}}{3}+(\alpha_k^2-\alpha_{k-1}^2)t\right)\right).
\end{align*}
Or, 
\begin{align*}
3s_k-2s_{k-1}-s_{k+1}=&\dfrac{3k-2k+2-k-1}{3}\left(-\alpha_1^2(T-\epsilon_1-\epsilon_2)+2\epsilon\alpha_1(T-\epsilon_2) \right. \\
& \left. \, -\epsilon^2(T-\epsilon_2)+L\right)\\
=&\dfrac{1}{3}\left(-\alpha_1^2(T-\epsilon_1-\epsilon_2)+2\epsilon\alpha_1(T-\epsilon_2)-\epsilon^2(T-\epsilon_2)+L\right),
\end{align*}
which is negative for $\alpha_1$ sufficiently large, as well as $\alpha_k^2-\alpha_{k-1}^2$ and $s_{k+1}-s_k$, while $\alpha_{k-1}-\alpha_k>0$. Therefore, $f_k(\xi_{k}^{-}(t),t)f_{k-1}^{-1}(\xi_k^{+}(t),t)$ tends to zero as $\alpha_1$ tends to infinity. 
Moreover, 
\begin{align*}
f_k(\xi_{k}^{-}(t),t)f_{k+1}(\xi_k^{-}(t),t)=&\exp\left(-\alpha_k\left(\frac{s_{k+1}-s_k}{3}\right)\right. \\
& \left. -\alpha_{k+1}\left(\frac{2s_k-2s_{k+1}}{3}\right)+(\alpha_k^2-\alpha_{k+1}^2)t \right) \\
=&\exp\left((\alpha_k-2\alpha_{k+1})\left(\frac{s_k-s_{k+1}}{3}\right) -\alpha_{k+1}(\alpha_k^2-\alpha_{k+1}^2)t \right),
\end{align*}
which tends to zero as $\alpha_1$ tends to infinity. 

Finally, let us consider the case where $k=N$. For the same reasons as in the case $2\leq k \leq N-1$, it is sufficient to consider the case $f_{N}(\xi_{k}^{-}(t),t)f_{N-1}^{-1}(\xi_k^{+}(t),t)$. We have, 
\begin{align*}
f_{N}(\xi_{k}^{-}(t),t)f_{N-1}^{-1}(\xi_k^{+}(t),t)=&\exp\left(-\alpha_N\left(\frac{s_N-s_{N-1}}{3}\right) \right. \\
& \left. \, \qquad +\alpha_{N-1}\left(4\frac{s_N-s_{N-1}}{3}+(\alpha_N^2-\alpha_{N-1}^2)t\right)\right)\\
=&\exp\left((4\alpha_{N-1}-\alpha_N)\left(\frac{s_N-s_{N-1}}{3}\right) \right. \\ 
& \left. \, \qquad  +\alpha_{N-1}(\alpha_N^2-\alpha_{N-1}^2)t\right),
\end{align*}
which tends to zero as $\alpha_1$ tends to infinity since $s_N-s_{N-1}<0$, $4\alpha_{N-1}-\alpha_N>0$, and $\alpha_N^2-\alpha_{N-1}^2<0$.

\endpf

We can now prove that $G$ satisfies $\Phi(x,T)\geq L, (x,t)\in (0,L)\times (T-\epsilon_2,T)$. We obtain the rough estimate on $G$, from (\ref{soln1}), (\ref{solnk}), and from the definition of the width of a soliton, 
\begin{equation*}
G(x,t)\geq \sum_{k=1}^N \dfrac{A_k(t)\alpha_k}{8} \mathds{1}(x)_{[\alpha_k^2t + \sigma_k(t)-w(\alpha_k)/2,\alpha_k^2t + \sigma_k(t)+w(\alpha_k)/2]},
\end{equation*}
where $w(\alpha)$ is the width defined by (\ref{width}). From Lemma \ref{lemAk}, we suppose that $\alpha_1$ is large enough so $A_k(t) \leq 1/2$, $\forall t\in [0,T]$ and $1\leq k \leq N$. We have,
\[
G(x,t)\geq \sum_{k=1}^N \dfrac{\alpha_k}{16} \mathds{1}(x)_{[\alpha_k^2t + \sigma_k(t)-w(\alpha_k)/2,\alpha_k^2t + \sigma_k(t)+w(\alpha_k)/2]},
\]
and, therefore,
\begin{align}
\dfrac{\partial \Phi}{\partial t}(x,t)=& G(\Phi(x,t),t) \nonumber \\
\geq& \sum_{k=1}^N \dfrac{\alpha_k}{16} \mathds{1}(\Phi(x,t))_{[\alpha_k^2t + \sigma_k(t)-w(\alpha_k)/2,\alpha_k^2t + \sigma_k(t)+w(\alpha_k)/2]}, \label{disp} \\
\geq& \sum_{k=1}^N \dfrac{\alpha_k}{16} \mathds{1}(x)_{[\alpha_k^2t + \sigma_k(t)-w(\alpha_k)/2,\alpha_k^2t + \sigma_k(t)+w(\alpha_k)/2]}, \label{nondisp}
\end{align}
the last line coming from the fact that, since all the characteristic functions are positive, the displacement of $\Phi$ will always be on the right. Therefore, the displacement of $\Phi$ will be greater in (\ref{disp}) than in (\ref{nondisp}), since the characteristic functions follow, even for a brief moment, the displacement of $\Phi$.

The displacement $\Phi$ under the flow of $G$ is then easily estimated from (\ref{nondisp}) 
\[
\ln\left(\sqrt{2\alpha_k}\left(1+\sqrt{1-\frac{1}{2\alpha_k}}\right)\right) / 4 \alpha_k^2,
\] 
that is, the height times the width divided by the speed $\alpha_k$ of each soliton. Consequently, each point $x\in [0,L]$ under the action of the flow of $G$ will move of at least $\ln\left(\sqrt{2\alpha_k}\left(1+\sqrt{1-\frac{1}{2\alpha_k}}\right)\right) / 4 \alpha_k^2$ to the right each time a soliton passes through the domain $[0,L]$. Since 
\[
 N=\ceil[\bigg]{4L \alpha_1^2 / \ln\left(\sqrt{2\alpha_1}\left(1+\sqrt{1-\frac{1}{2\alpha_1}}\right)\right)},
\]
 we obtain that $\Phi(x,T)\geq L, (x,t)\in (0,L)\times (T-\epsilon_2,T)$. 

\endpf

\textbf{Acknowledgements}

The author would like to thank Jean-Michel Coron for the helpful discussions and for the introduction to this problem.

\bibliographystyle{plain}
\bibliography{biblio}

\end{document}